\documentclass[11pt,a4paper,psamsfonts]{amsart}

\usepackage{amsmath}
\usepackage{fancyhdr}
\usepackage{appendix}
\usepackage{amssymb,amscd,amsxtra,calc}
\usepackage{mathrsfs}
\usepackage{cmmib57}
\usepackage{multirow}
\usepackage[all]{xy}
\usepackage{longtable}
\usepackage{enumitem}
\usepackage{mathtools}
\usepackage{comment}

\usepackage[colorlinks,linkcolor=blue,anchorcolor=blue,citecolor=green,backref=page]{hyperref}

\theoremstyle{plain}
    \newtheorem{thm}{Theorem}[section]
    
    \newtheorem{claim}[thm]{Claim}
    \newtheorem{conjecture}[thm]{Conjecture}
    \newtheorem{corollary}[thm]{Corollary}
    
    \newtheorem{lemma}[thm]{Lemma}
    \newtheorem{proposition}[thm]{Proposition}
    
    \newtheorem{theorem}[thm]{Theorem}

\theoremstyle{definition}

    \newtheorem*{notation*}{Notation and Terminology}
    \newtheorem{remark}[thm]{Remark}
\theoremstyle{remark}

    \newtheorem{setup}[thm]{}

\newcommand{\C}{\mathbb{C}}

\newcommand{\PP}{\mathbb{P}}

\newcommand{\Q}{\mathbb{Q}}

\newcommand{\Aut}{\operatorname{Aut}}

\newcommand{\Gal}{\operatorname{Gal}}

\newcommand{\PGL}{\operatorname{PGL}}

\newcommand{\Sing}{\operatorname{Sing}}

\newcommand{\Pic}{\operatorname{Pic}}

\newcommand{\reg}{\mathrm{reg}}

\usepackage[colorinlistoftodos]{todonotes}

\newcommand{\nc}{\newcommand}
\nc{\cH}{{\mathcal H}}
\nc{\cA}{{\mathcal A}}
\nc{\cG}{{\mathcal G}}
\nc{\cC}{{\mathcal C}}
\nc{\cO}{{\mathcal O}}
\nc{\cI}{{\mathcal I}}
\nc{\cB}{{\mathcal B}}
\nc{\cY}{{\mathcal Y}}
\nc{\cK}{{\mathcal K}}
\nc{\cX}{{\mathcal X}}
\nc{\cS}{{\mathcal S}}
\nc{\cE}{{\mathcal E}}
\nc{\cF}{{\mathcal F}}
\nc{\cZ}{{\mathcal Z}}
\nc{\cQ}{{\mathcal Q}}
\nc{\cN}{{\mathcal N}}
\nc{\cP}{{\mathcal P}}
\nc{\cL}{{\mathcal L}}
\nc{\cM}{{\mathcal M}}
\nc{\cT}{{\mathcal T}}
\nc{\cW}{{\mathcal W}}
\nc{\cU}{{\mathcal U}}
\nc{\cJ}{{\mathcal J}}
\nc{\cV}{{\mathcal V}}

\nc{\bH}{{\mathbb H}}
\nc{\bA}{{\mathbb A}}
\nc{\bG}{{\mathbb G}}
\nc{\bC}{{\mathbb C}}
\nc{\bO}{{\mathbb O}}
\nc{\bI}{{\mathbb I}}
\nc{\bB}{{\mathbb B}}
\nc{\bY}{{\mathbb Y}}
\nc{\bK}{{\mathbb K}}
\nc{\bX}{{\mathbb X}}
\nc{\bS}{{\mathbb S}}
\nc{\bE}{{\mathbb E}}
\nc{\bF}{{\mathbb F}}
\nc{\bZ}{{\mathbb Z}}
\nc{\bQ}{{\mathbb Q}}
\nc{\bN}{{\mathbb N}}
\nc{\bP}{{\mathbb P}}
\nc{\bL}{{\mathbb L}}
\nc{\bM}{{\mathbb M}}
\nc{\bT}{{\mathbb T}}
\nc{\bW}{{\mathbb W}}
\nc{\bU}{{\mathbb U}}
\nc{\bD}{{\mathbb D}}
\nc{\bJ}{{\mathbb J}}
\nc{\bV}{{\mathbb V}}
\nc{\bbZ}{{\mathbb Z}}
\nc{\bR}{{\mathbb R}}
\nc{\fr}{{\rightarrow}}
\nc{\co}{{\nabla}}

\nc{\cu}{{\overlineline{\nabla}}}

\title [morphisms]{Morphisms from a very general hypersurface}

\author{Yongnam Lee, Yujie Luo and De-Qi Zhang}
\date{}

\address{Center for Complex Geometry, Institute for Basic Science (IBS), 55 Expo-ro, Yuseong-gu, Daejeon 34126, Korea, and
\newline \hspace*{3mm} Department of Mathematical Sciences, KAIST, 291 Daehak-ro, Yuseong-gu, Daejeon 34141, Korea}

\email{ynlee@kaist.ac.kr}

\address{Department of Mathematics, National University of Singapore, Singapore 119076, Republic of Singapore}

\email{yujieluo96@gmail.com~and~lyj96@nus.edu.sg}

\address{Department of Mathematics, National University of Singapore, Singapore 119076, Republic of Singapore}

\email{matzdq@nus.edu.sg}


\begin{document}

\begin{abstract}
Let $X$ be a very general hypersurface of degree $d$ in the projective $(n+1)$-space with $n \ge 3$, and $f: X \to Y$ a non-birational surjective morphism to a normal projective variety $Y$. We first prove that $Y$ is a klt Fano variety 
if $\deg f \ge C$ for some constant $C = C(n, d)$ depending only on $n$ and $d$. Next we prove an optimal upper bound $\deg f \le \deg X$ provided that $Y$ is factorial, $\deg f$ is prime and $\deg f \ge E(n)$ for some constant $E(n)$ (with $E(n) = n(n+1)$ when $Y$ is smooth). As a corollary, we show that $Y\cong\bP^n$ under some conditions on $Y$ and $\deg f$. 
\end{abstract}

\maketitle

\begin{center}
Dedicated to James M\textsuperscript{c}Kernan on the occasion of his 60th birthday.
    
\end{center}

\tableofcontents


%
%
%
%
\section{Introduction}

We work over the field $\C$ of complex numbers.

To motivate our result, we recall that if $X$ is a smooth projective variety of general type, then the set of dominant rational maps $X\dashrightarrow Y$ from $X$ to a smooth variety $Y$ of general type, up to birational equivalence of $Y$, is a finite set. The proof was given by Maehara \cite{Ma} under the assumption of effective birationality of pluri-canonical maps of varieties of general type. This effective birationality was proved later (cf.~Hacon and M\textsuperscript{c}Kernan \cite[Theorem 1.1, Corollary 1.4]{HM}, Takayama \cite {Ta}, and Tsuji \cite{Ts}).

If we allow the co-domain $Y$ to be uniruled (or even rational), then the above finiteness result is not true since
every variety has
a generically finite map to a projective space while the latter has self maps of unbounded degrees.
In this paper, we assume that the domain $X$ is a very general hypersurface.

In connection with the {\bf famous conjecture}:

\par \noindent
`{\it Let $Y$ be a Fano manifold of Picard number 1. Suppose that $Y$ admits a non-isomorphic surjective endomorphism. Then $Y$ is a projective space.}' 

\par \noindent
which has been proved in dimension $3$ by Amerik \cite{Am97} and Amerik-Rovinsky-Van de Ven \cite{ARV}, independently by Hwang-Mok \cite{HM03}, and even in arbitrary characteristics, coprime with $\deg f$, by Kawakami-Totaro \cite{KT}, it is natural to
consider the following about morphisms from a very general hypersurface (cf.~\cite[Conjecture B]{Am97}):

\begin{conjecture}\label{conj: 2}
Let $X \subset \bP^{n+1}$ be a very general hypersurface with $n \ge 2$ and $\Pic(X)=\bZ$ (this holds when $n \ge 3$ or $\deg X \ge 4$), and $f: X \to Y$ a non-birational surjective morphism to a smooth projective variety $Y$ of positive dimension. 
Then there is a constant $C(X)$ depending on discrete invariants of $X$, such that if $\mathrm{deg}~f\ge C(X)$, then $Y$ is isomorphic to $\bP^n$.
\end{conjecture}

Conjecture \ref{conj: 2} in dimension $\le 3$ has been successfully proved by Amerik, and more recently, by Kawakami and Totaro. Indeed, see Amerik \cite[Theorems 1.1 and 1.2]{Am} when $n = 2$, or when $n=3$ and $Y$ is neither the Mukai–Umemura variety $V_{22}^s$ nor of type $V_5$; see \cite[Proposition 2.2]{Am97} for the last two Fano $3$-fold cases. See also \cite[Proof of Proposition 2.7, Proposition 3.8]{KT} (which works in arbitrary characteristics) and \cite[Corollary 2.2, Remark 2.4, Lemma 3.0]{Am}. We refer to Shao-Zhong \cite{SZ} for their results for morphisms between Fano manifolds, and more references and progress towards the famous conjecture and Conjecture \ref{conj: 2}.  

Here, we achieve a partial result for Conjecture \ref{conj: 2}. Our first main result (Theorem \ref{thm: main2}) is to show that $Y$ is a klt Fano variety if $\deg f$ is large enough.

\begin{theorem}\label{thm: main2} Let $X \subset \bP^{n+1}$ be a very general hypersurface with $n\ge 3$, and $f: X \to Y$ a non-birational surjective morphism to a normal projective variety $Y$. Then we have:
\begin{enumerate}
\item[(1)]
$Y$ has only klt singularities; either the canonical divisor $K_Y$ is ample, or $-K_Y$ is ample.
\item[(2)]
There exists a constant $C = C(n, d)$ such that if $\mathrm{deg}~f\geq C$, then $Y$ is a klt Fano variety.
\end{enumerate}
\end{theorem}

As remarked earlier, Conjecture \ref{conj: 2} in dimension $3$ has been confirmed by Amerik and more recently by Kawakami and Totaro, both thanks to the known classification list of smooth Fano $3$-folds due to Iskovskikh and Mori-Mukai. However, a classification of Fano manifolds in higher dimensions is unlikely to be feasible. Therefore, for the morphism $f: X \to Y$ in the above theorem, 
assuming further that $\deg f$ is large with additional conditions, we prove in Theorem \ref{thm: e lc-1} and Corollaries \ref{cor: p3} and \ref{cor: pn} that $Y \cong \bP^n$, thus confirming Conjecture \ref{conj: 2}. This also gives the precise constant $C(n, d)$ and even removes or weakens the smoothness assumption on $Y$. 

Before stating our next main results, we introduce some notation. Let $X \subset \bP^{n+1}$ be a very general hypersurface of degree $d$ with $n\ge 2$, and $f: X \to Y$ a non-isomorphic surjective morphism of degree $m$ to a normal projective variety $Y$.  Consider the restriction map $f_{H_X} := f|_{H_X}$ where $H_X: =H\cap X$ with $H \subset \bP^{n+1}$ a hyperplane, it follows that $f_{H_X}$ is finite. Let $H_Y: =f(H_X)$. Then, by the projection formula, $f_* f^*H_Y= (\deg f)H_Y$. Therefore, $\deg f_{H_X}$ divides $\deg f$ when $\mathrm{Pic}(X)=\mathbb{Z}$; indeed, writing $f^*H_Y \sim rH_X$ for some positive integer $r$ and using the projection formula to push forward by $f_*$, we then have $\deg f = r \deg (f_{H_X}: H_X \to H_Y)$.
We define a linear subspace $S(f)$ of $|H_X|$ below (which is the pullback of a complete linear system on $Y$; see Proposition \ref{prop: Weil gen})
$$S(f): =\{ H'\in |H_X| \, ; \, \deg f|_{H'}=\deg f\}.$$

The $Y$ in Theorem \ref{thm: e lc-1} below, which is the image of smooth $X$, has at worst klt singularities (cf.~Proposition \ref{prop: basic}). So, it is quite natural to further assume that $Y$ has at worst $\varepsilon$-lc singularities (cf.~Birkar \cite[2.7]{Bir} for its definition). Recall that a canonical singularity is $\varepsilon$-lc with $\varepsilon = 1$. The $Y$ in Theorem \ref{thm: e lc-1} or \ref{thm: e lc-2} below, being klt, has at worst canonical singularities if $Y$ is further factorial. Let $Y_{\rm reg} := Y \setminus \Sing Y$ be the smooth locus of $Y$.

\begin{theorem}\label{thm: e lc-1}
Let $X \subset \bP^{n+1}$ be a very general hypersurface of degree $d$
and $f: X \to Y$ a surjective morphism to a normal projective variety $Y$ with at  worst $\varepsilon$-lc singularities. Assume either $n \ge 3$ and $d \ge 3$, or $n = 2$ and $d \ge 4$.  Assume further the fundamental group $\pi_1(Y_{\rm reg})=(1)$ (which holds if $Y$ is smooth or $\deg f$ is prime by Proposition \ref{prop: basic}).
Then we have:
\begin{enumerate}
\item 
$\deg f\le D_1 d$ for some constant $D_1=D_1(n, \varepsilon)$, if $-K_Y$ is ample and $S(f)\ne\emptyset$.
\item $\deg f\le {\rm max}\{D_1 d, D_2\}$ for some constant $D_2=D_2(n, \varepsilon, d)$, if $\deg f$ is prime.
\item If $\deg f \ge d$, then $\deg f = d$ and $Y = \bP^n$, unless ${\rm Codim}~S(f)\ge 2$ in $|H_X|$.
\end{enumerate}
\end{theorem}

Theorem \ref{thm: e lc-2} below says that in Theorem \ref{thm: e lc-1}, if we assume $Y$ is factorial then $D_1=1$, and $D_2$ depends only on $n$. The upper bound in Theorem \ref{thm: e lc-2} below is optimal because of the projection map (see Remark \ref{rem: on main thm}).

An upper bound of $\deg f$ without a condition like $\deg f$ being prime or so in Theorem \ref{thm: e lc-1} or \ref{thm: e lc-2}, seems unlikely. Indeed, a projection $f_0: X \to \bP^n$ composed with quotient maps $f_s: \bP^n \to Y_s := \bP^n/G_s$ for a finite subgroup $G_s \subset \PGL_{n+1}(\bC)$, will give morphisms from $X$ with unbounded degrees, with $Y_s$ being smooth (and hence is isomorphic to $\bP^n$) when $G_s$ is generated by pseudo-reflections, or being singular otherwise.

\begin{theorem}\label{thm: e lc-2}
Let $X \subset \bP^{n+1}$ be a very general hypersurface of degree $d$
and $f: X \to Y$ a surjective morphism 
to a factorial normal projective variety $Y$. Assume $n \ge 3$, $d \ge 3$, and $\pi_1(Y_{\rm reg})=(1)$ (which holds if $Y$ is smooth or $\deg f$ is prime by Proposition \ref{prop: basic}).
Then we have:
\begin{enumerate}
\item 
$\deg f \le d$, if $-K_Y$ is ample and $S(f)\ne\emptyset$.
\item $\deg f \, \le \, {\rm max}\{d, \,\,  D_2/d^{1/(n-1)}\}$ for some constant $D_2 = D_2(n)$,
if $\deg f$ is prime.
\item 
$\deg f \le d$, if $\deg f$ is prime and $d \ge D_2^{(n-1)/n}$.
\end{enumerate}
\end{theorem}

As a consequence of Theorems~\ref{thm: e lc-1} and \ref{thm: e lc-2} (or Theorem \ref{thm: e lc'}), we show that $Y\cong\bP^n$ under some conditions on $Y$ and $\deg f$ (see Corollaries~\ref{cor: p3} and \ref{cor: pn}).

\begin{remark}\label{rem: on main thm}
$ $
\begin{enumerate}
\item
If $f: X \to Y = \bP^n$ is a projection morphism from a point in $\bP^{n+1}\setminus X$,
then $\deg f = d = \deg X$ and ${\rm Codim}~S(f) = 1$ in $|H_X|$. Indeed, the pullbacks to $X$ of hyperplanes in $Y$ give $S(f)$. In particular, the upper bound in Theorem \ref{thm: e lc-2} is optimal. On the other hand, we have $\deg f \ge d-n$, by Lemma \ref{lem: degree}, when $X$ is of general type (i.e., $d \ge n+3$).

\item 
In Theorem \ref{thm: e lc-2}, we can take $D_2 = (D_2')^{1/(n-1)}$ in terms of upper bound $D_2'$ of $(-K_Y)^n$ (see Theorem \ref{thm: e lc'}). So the condition in Theorem \ref{thm: e lc-2} (3) is $d \ge D_2^{(n-1)/n} = (D_2')^{1/n}$ ($ = n(n+1)$ when $Y$ is smooth, see \cite[Corollary 1]{KMM}).

\item 
In Theorem \ref{thm: e lc-1}, when $\deg f = d$, we expect that $f$ is a projection morphism, and hence the exceptional case in Theorem \ref{thm: e lc-1} (3) is expected not to occur.

\item 
In Theorem \ref{thm: e lc-1}, when $\deg f \ge d$, $``{\rm Codim} \, S(f) \le 1$ in $|H_X|$'' holds if and only if $\deg f = d$, $Y = \bP^n$ and $\cO(f(H_X)) = \cO(1)$ (cf. Proposition \ref{prop: Weil gen}).

\item 
In Theorem \ref{thm: e lc-2}, the condition $n \ge 3$ and $d \ge 3$ can be weakened as either $n \ge 3$ and $d \ge 3$, or $n = 2$ and $d \ge 4$ (just the same as in Theorem \ref{thm: e lc-1}).

\item
From the Lefschetz theory and the irreducibility of the monodromy action on vanishing cycles (cf.~Voisin \cite[Chapters 1-3]{Vo}), the non-birational curve image of a very general plane curve of degree $d\ne 3$ is a rational curve. Hence we consider only the case $n \ge 2$ in our theorems. 

\end{enumerate}
\end{remark}

Our proofs of the theorems combines Hodge theory and some birational geometry including Hacon-M\textsuperscript{c}Kernan-Xu's work \cite{HMX} and Birkar's work \cite{Bir}.
Fujita's theory (\cite{Fuj82}, \cite{Fuj89}) on (quasi-) polarized projective varieties and $\Delta$-genus is also crucial in proving Proposition \ref{Prop: HP} (essentially from H\"oring-Peternell \cite[Theorem 3.2]{HP}) and Proposition \ref{prop: Weil gen}, which in turn, are needed for proving Theorems~\ref{thm: e lc'}, \ref{thm: e lc-1} and \ref{thm: e lc-2}; see also Kobayashi-Ochiai \cite{KO}.

The method presented here can be used to obtain similar results for a very general hypersurface of large degree in a Fano manifold of Picard number one; see Remarks \ref{rem: basic-fano} and \ref{rem: fano main thm}.

It is not clear at the moment how to treat the case where $f$ is not a morphism.

\vspace{2mm}

{\bf Acknowledgements.}
The authors would like to thank IBS and KIAS in Korea, NUS in Singapore, and Simons Center for Geometry and Physics in the United States for their warm hospitality, and E. Amerik, G. Zhong and the referee for the valuable suggestions. Y. Lee is supported by the Institute for Basic Science IBS-R032-D1, Y. Luo by a PTA fellowship of NUS, and D. -Q. Zhang by ARF: A-8002487-00-00 of NUS.
\medskip

\section{Preliminary results}

\begin{setup} {\bf Terminology and Notation}
We adopt the standard terminology and notation as in \cite{Har}, and push-forward of Cartier divisors for a proper morphism (which preserves linear equivalence) as in \cite[Proposition 1.4]{Ful}. 
We also refer to \cite{KM} or \cite[2.7]{Bir} for the definitions of {\it Kawamata log terminal} (klt), {\it canonical} and {\it $\varepsilon$-lc} (with $\varepsilon > 0$) singularities.

Let $f: X\dashrightarrow Y$ be a dominant rational map between two varieties of the same dimension. The {\it degree} of $f$, denoted as $\mathrm{deg}~f$, is defined as the degree of the function field extension $K(X)/K(Y)$. It is also equal to the cardinality of $f^{-1}(y)$ for a general point $y \in Y$.

For a finite surjective morphism $f: X \to Y$ between normal projective varieties, we use very often the {\it ramification divisor formula} 
$$K_X = f^*K_Y + R_f$$ 
and pullback $f^*D$ for a Weil divisor $D$ on $Y$. Precisely, we restrict $D$ on the smooth locus of $Y$, use $f^*$ to pull the restriction back, and then take the Zariski closure. This pullback definition coincides with the usual definition when $D$ is a $\Q$-Cartier divisor.

A hypersurface $X \subset \bP^{n+1}$ of degree $d$ is {\it very general} if it is outside the union of countably many proper closed subsets of hypersurfaces of degree $d$.
\end{setup}

We begin with the following key proposition.

\begin{proposition}
\rm ({\bf Key Proposition}) \label{prop: basic}
Let $X \subset \bP^{n+1}$ be a very general hypersurface of degree $d$, $Y$ a normal projective variety of positive dimension, and $f: X \to Y$ a surjective morphism. Suppose either $n \ge 3$ and $d \ge 3$, or $n = 2$ and $d \ge 4$.
Then the following hold
(indeed, (1) - (3) hold whenever either $n \ge 3$ and $d \ge 1$, or $n = 2$ and $d \ge 4$).
\begin{enumerate}
\item 
$\Aut(X) = (1)$.
\item  
Let $H_X$ be a general section of $\cO_{\bP^{n+1}}(1)$ restricted to $X$. Then $\Pic X = \bZ[H_X]$ 
(so the Picard number $\rho(Y) = 1$ and every member of $|H_X|$ is reduced and irreducible), and the fundamental group $\pi_1(X) = (1)$ (so there is no non-trivial torsion line bundle on $X$).
\item $f$ is a finite morphism,
\item $Y$ is $\bQ$-factorial and klt.
\end{enumerate}
Suppose further that $f$ is non-birational. Then:
\begin{enumerate}[start=5]
\item 
If $d \le n+2$ then $Y$ is a klt Fano variety.
\item 
$H^i(Z, \cO_Z)=0$ for any smooth model $Z$ of $Y$ and all $i >0$.
\item $p_g(Z) = 0$, and $\chi(\cO_Z)=1$. 
\item
$\pi_1(Z) = (1) = \pi_1(Y)$.
\item
For $Y_{\reg} := Y \setminus \Sing Y$, we have $|\pi_1(Y_{\reg})|$ being a finite number dividing $\deg f$.
\item 
Suppose $\deg f$ is prime. Then $\pi_1(Y_{\reg}) = (1)$. So for Weil integral divisors $M_j$ on $Y$ with $sM_1 \sim sM_2$ for some integer $s \ge 1$, we have $M_1 \sim M_2$.
\end{enumerate}
\end{proposition}

\begin{proof}
By the assumption on $n$ and $d$, we have $\Aut(X) = (1)$ by \cite[page 347]{MM}, and also
$$\pi_1(X) = \pi_1(\bP^{n+1}) = (1)$$ 
and the restriction map 
$$\Pic(\bP^{n+1}) \to \Pic(X)$$ is an isomorphism by the Lefschetz hyperplane theorem (cf.~\cite{Noe} and \cite{Lef}). Then (2) follows. 

As in (2), write 
$$\bZ[H_X] = \Pic(X) = \Pic(X)/\Pic^0(X) = \mathrm{NS}(X)$$
where $H_X$ is a section of $\cO_{\bP^{n+1}}(1)$ restricted to $X$ and it is very ample.
Since $\dim\, Y > 0$, $\rho(Y) = 1 = \rho(X)$, hence $f^*$ induces an isomorphism $$\mathrm{NS}(Y)_{\bQ} \to \mathrm{NS}(X)_{\bQ}$$
of $\bQ$-coefficient Neron-Severi spaces, so $H_X$ is equal to some $f^*D_Y$. Thus (3) follows.

Since $f: X\to Y$ is a finite surjective morphism and $X$ is smooth (and hence $\bQ$-factorial and klt), we obtain (4) by \cite[Lemma~5.16 and Proposition~5.20]{KM}.

For (9), since $\dim f^{-1}(\Sing Y)$ is of codimension $\ge 2$ in the smooth variety $X$ and by \cite[Proposition 1.3]{cam91}, the image of $(1) = \pi_1(X) = \pi_1(f^{-1} \, Y_{\reg})$ in $\pi_1(Y_{\reg})$ is of index $\le \deg f$, so $|\pi_1(Y_{\reg})|$ is finite, and further, $X \to Y$ factors as $X \to Y' \to Y$, where $Y' \to Y$ is finite and quasi-\'etale and is the normalisation of $Y$ in the function field of the universal cover of $Y_{\reg}$.
Thus $|\pi_1(Y_{\reg})| = \deg(Y'/Y)$ divides $\deg f$, finishing the proof of  (9).

For (10), as in (9), we have $X \to Y' \to Y$ with $|\pi_1(Y_{\reg})| = \deg(Y'/Y)$ dividing $\deg f$. If $|\pi_1(Y_{\reg})| \ge 2$ then the primality of $\deg f$ implies that $X = Y'$, and hence $\Gal(Y'/Y) \le \Aut(X) = (1)$, which is absurd. Thus $\pi_1(Y_{\reg}) = (1)$.
If $sL \sim 0$ for some integer $s \ge 2$
while $L \not\sim 0$, taking a prime factor $q$ of $s$ and replacing $L$ by a multiple, we may assume that $qL \sim 0$ while $L \not\sim 0$. This relation induces a finite Galois $\bZ/q\bZ$-cover $Y' \to Y$, which is quasi-\'etale, i.e., \'etale over $Y_{\reg}$ by the purity of branch loci.
Thus $|\pi_1(Y_{\reg})|\ge q$, contradicting its triviality which we have just proved. This proves (10).

Next we prove (5) - (8).
By the ramification divisor formula,
$$K_X = f^*K_Y + R_f$$ with $R_f$ the ramification divisor.
Since $\rho(X) = 1$ and $\pi_1(X) = (1)$, we have linear equivalence
$$f^*K_Y \sim k_y H_X$$
for some integer $k_y$.

\begin{claim}\label{claim: Rf=0 and kY=0}
It is impossible that $k_Y = 0$ and $R_f = 0$.
\end{claim}

We prove the claim. Suppose the contrary that $k_y = 0$ and $R_f = 0$. Then $K_X \sim 0$, and $f$ is quasi-\'etale. Let 
$$X' \to X \to Y$$ be the Galois closure of $f$. 
Then $X' \to X$, like $f$, is still quasi-\'etale,
and hence \'etale since $X$ is smooth and by the purity of branch loci. Thus $X' = X$ because $\pi_1(X) = (1)$.
So $f: X \to Y$ is Galois with Galois group $G$.
Now $G \le \Aut(X) = (1)$. Hence $\deg f = |G| = 1$, contradicting the non-birationality assumption of $f$. This proves Claim \ref{claim: Rf=0 and kY=0}.

We return back to the proof of the proposition.

If $k_y < 0$, then $-K_Y$ is ample. Thus $Y$ is a klt Fano variety, and hence both $Y$ and its smooth model $Z$ are rationally connected, so (6) - (7) clearly hold in this case. Since smooth and hence normal rationally connected varieties are simply connected \cite[Theorem 3.5]{cam91}, the assertion (8) also holds in this case.

Thus, it remains to show (5) - (8) under the assumption that either $k_y > 0$,  or $k_y=0$ and $R_f \ne 0$. In either case, from the ramification divisor formula we see that $$K_X\sim (d-n-2)H_X$$ is a nonzero effective divisor, and hence ample since $\rho(X) = 1$.
Thus $d \ge n+3 \ge 5$ (this also proves (5)). 
Further, $K_X$ is very ample and $X$ is of general type.
Replacing $Y$ by its smooth model $Z$ and $X$ with a higher smooth birational model $W$, our $f$ induces a generically finite surjective morphism 
$$h:W\to Z.$$ Since $X$ is simply connected, so is $W$, hence $|\pi_1(Z)| = |\pi_1(Z): h(\pi_1(W))|$ is bounded by the number of connected components of a general fibre, i.e., by $\deg f$  (cf. \cite[Proposition 1.3]{cam91}). Let 
$$Z'\to Z$$ be the universal cover. Then $W\to Z$ factors as 
$$W \to Z' \to Z$$ because $W$ is simply connected. 

Since $d \ge n+3$, by \cite[Proposition~2.1 ]{Am} or \cite[Remark~3.3, Example 1.2]{BDet}, either $p_g(Z) = 0$, or $p_g(Z) = p_g(W)$ ($ = p_g(X)$).

\begin{claim}\label{claim: geometric genus Z is W}
It is impossible that $p_g(Z) = p_g(W)$ ($ = p_g(X)$).
\end{claim}

Suppose to the contrary that $p_g(Z) = p_g(W)$. Now the natural injection, via $h^*$,
$$H^0(Z, K_{Z}) \to H^0(W, K_W) = H^0(X, K_X)$$ is an isomorphism, so the canonical rational map $\Phi_{|K_W|}$ 
(which is birational, since $K_X$ is very ample)
factors as the composition of rational maps
$$W \to Z \overset{\Phi_{|K_Z|}}\dashrightarrow \Phi_{|K_Z|}(Z),$$
or equivalently as the composition of rational maps
$$W \rightarrow X \overset{f}\to Y \dashrightarrow Z \dashrightarrow \Phi_{|K_Z|}(Z)$$
where $W \to X$ is birational and $Y \dasharrow Z$ is the inverse of the birational map $Z \dasharrow Y$.
Hence $\deg f = 1$, a contradiction. This proves Claim \ref{claim: geometric genus Z is W}.

We return back to the proof of the proposition. By Claim \ref{claim: geometric genus Z is W}, $\dim H^n(Z, \cO_Z) = p_g(Z) = 0$. Applying the Kodaira vanishing to the cohomologies of the exact sequence 
$$0\to \cO_{\bP^{n+1}}(-X)\to \cO_{\bP^{n+1}}\to \cO_X\to 0$$ 
we have $H^i(X, \cO_X)=0$ for all $i=1, \dots, n-1$, so the same holds for $Z'$ and $Z$, since the latter ones are dominated by $X$. Combining the above all, we get $\chi(\cO_Z) = 1$. This proves (6) - (7) in this case too.

\begin{claim}\label{claim: geometric genus non zero}
$p_g(Z') =0$.  
\end{claim}

We prove the claim. Suppose the contrary that $p_g(Z') \ge 1$. Then $\deg(Z'/Z) \ge 2$. Applying \cite[Proposition~2.1 ]{Am} or \cite[Remark~3.3, Example 1.2]{BDet} and arguing as above, we have 
$$p_g(Z') = p_g(W) = p_g(X)$$ and $W \to Z'$ is a birational morphism.
Now $G = \Gal(Z'/Z)$ acts birationally on $X$ and hence
biregularly on the canonical model of $X$, which is $X$ itself. Thus $G \le \Aut(X) = (1)$. Hence $Z' = Z$, a contradiction. This proves Claim \ref{claim: geometric genus non zero}

We return back to the proof of part (8) of the proposition. By Claim \ref{claim: geometric genus non zero}, $p_g(Z') = 0$. So $\chi(\cO_{Z'}) = 1$.
Since $Z' \to Z$ is \'etale, $\deg(Z'/Z)| \, \chi(\mathcal{O}_{Z})| = |\chi(\cO_{Z'})|$ ($= 1$). Hence $\deg(Z'/Z) = 1$, i.e., $Z' = Z$.
Thus $Z$ is simply connected.
The same holds for $Y$ since we may assume $Z \to Y$ is a surjective birational morphism and $Y$ is normal
(cf.~\cite[Proposition 1.3]{cam91}). This proves (8) in this case too.
\end{proof}

\begin{remark}\label{rem: basic-fano}
All assertions in Proposition~\ref{prop: basic}, except possibly (4), (5), (8) and (10), can be generalized from $\bP^{n+1}$ to any nondegenerate Fano manifold $P$ of $\dim\, P\ge 3$ with the Picard number $\rho(P)=1$, if $X$ is a very general hypersurface section of $P$ and $X$ is of general type. If we further assume that $\Aut(X)= (1)$, then all the assertions in Proposition~\ref{prop: basic} can be generalized except (5).

Indeed, let $f: X\to Y$ be a non-birational surjective morphism to a normal projective variety $Y$ of positive dimension. Assuming that $X$ is of general type and $P$ is a Fano manifold with $\rho(P)=1$, we can apply again the irreducibility of the monodromy action on vanishing cycles (cf. \cite[Chapters 1-3]{Vo}) and Kodaira vanishing to show $H^i(\cO_Z)=0$ for any smooth model $Z$ of $Y$ and all $i>0$. Other proofs are basically the same.
\end{remark}

\section{Proof of Theorem~\ref{thm: main2}}\label{proof of main2}

\begin{proof}[{\bf Proof of Theorem~\ref{thm: main2}}]
We may assume that $\mathrm{dim}~Y>0$. By Proposition~\ref{prop: basic}, $f$ is finite (so $m := \deg\, f \ge 2$), and $Y$ is $\bQ$-factorial and klt. 
We use the notation in Proposition~\ref{prop: basic}:
$${\rm Pic}\, X = \bZ[H_X],\ \ \ f^*K_Y \sim k_y H_X.$$
If $k_y < 0$, then $Y$ is a klt Fano variety, and the theorem holds.
Thus, we may assume $k_y \ge 0$.

\begin{claim}\label{claim: ample canonical divisor}
$K_Y$ is a $\bQ$-Cartier ample divisor.
\end{claim}

The above claim would complete the first part of the theorem. Now we prove this claim.

If $k_y > 0$, then $f^*K_Y$ and hence $K_Y$ are ample $\Q$-Cartier divisors. Thus we may assume $k_y = 0$. By the projection formula,
$$mK_Y = f_*f^*K_Y \sim f_*(k_y H_X) = 0.$$
Let $Y' \to Y$ be the index-1 cover which is quasi-\'etale, so that $K_{Y'} \sim 0$ and $Y'$ has only canonical singularities. As in the proof of Proposition~\ref{prop: basic}, $\pi_1(X) = (1)$ implies that $X \to Y$ factors as 
$$X \to Y' \to Y.$$
Let $Y''$ be a smooth model of $Y'$.
Then the Kodaira dimension $\kappa(X) \ge \kappa(Y'') = 0$. Hence $K_X \sim (d-n-2)H_X$ is effective, so $d \ge n+2 \ge 5$.

Since $p_g(Y'') = 1 > 0$, the morphism $X \to Y'$ is birational, or else contradicting Proposition~\ref{prop: basic} applied to $X \to Y'$.
Hence $\kappa(X) = \kappa(Y'') = 0$.
Thus $d = n+2$ and $K_X \sim 0$.
So $f: X \to Y$ is quasi-\'etale.
As in Proposition~\ref{prop: basic}, $\pi_1(X) = (1)$ implies that the Galois closure of $f$ is $f$ itself. Then $\Gal(X/Y) \le \Aut(X) = (1)$ by Proposition~\ref{prop: basic}, which implies that $\deg f = 1$, a contradiction. This finishes the proof of Claim~\ref{claim: ample canonical divisor}.

We return back to the proof of the theorem.
Now $K_Y$ is ample. Since $\rho(X) = 1$,
$$(d-n-2) H_X \sim K_X = f^*K_Y + R_f$$ is ample too. So $d \ge n+3$
and $X$ is of general type.

By Hacon-M\textsuperscript{c}Kernan-Xu \cite[Theorem 1.3]{HMX} (or Birkar-Zhang \cite[Theorem 1.3]{BZ}), there is a constant $C_1 = C_1(n)$ such that
$C_1 K_Y$ is birational, so 
$C_1K_Y > J_Y$
for some nef and big pullback $J_Y$ of an ample Cartier divisor via some pluri-canonical morphism of $Y$. Thus $C_1K_X >  f^*J_Y$. Since $K_X \sim (d-n-2) H_X$ and $H_X^n = d$, if we set $$C := C(n, d) = (C_1 K_X)^n$$
we have $$C = d C_1^n (d-n-2)^n.$$
Now $C_1K_X >  f^*J_Y$ and the claim below imply
$$C = (C_1 K_X)^n \ge  (\deg\, f) J_Y^n \ge \deg\, f .$$

\begin{claim}
Let $V$ be a projective variety of dimension $n \ge 1$, with $\Q$-Cartier nef divisors $L_1, L_2$. Suppose $M := L_1 - L_2$ is pseudo-effective. Then 
$$L_1^n \ge L_2^n.$$
\end{claim}

We prove the claim.
The nefness of $L_j$ implies
$$L_1^n = L_1^{n-1} (L_2 + M) \ge L_1^{n-1} L_2 =
L_1^{n-2} (L_2 + M) L_2 \ge L_1^{n-2} L_2^2 \ge \cdots \ge L_2^n.$$
This proves the claim and also the theorem.
\end{proof}

\begin{remark}\label{rem: fano main thm}
Theorem~\ref{thm: main2} can be generalized from $\bP^{n+1}$ to any nondegenerate Fano manifold $P$ of $\dim\, P\ge 4$ with the Picard number $\rho(P)=1$, if $X$ is a very general hypersurface section of $P$ and $X$ is of general type. Indeed, let $f: X \to Y$ be a non-birational surjective morphism to a normal projective variety $Y$ of positive dimension. Then there is a constant $C = C(n, d)$ where $d$ is the degree of the hypersurface section $X$ such that if $\deg\, f\ge C$ then $Y$ is a klt Fano variety. The proof is the same as above by Remark~\ref{rem: basic-fano} and by Hacon-M\textsuperscript{c}Kernan-Xu \cite[Theorem 1.3]{HMX} (or Birkar-Zhang \cite[Theorem 1.3]{BZ}).
\end{remark}

\begin{remark}
Assume $Y$ is smooth in Theorem~\ref{thm: main2}. If $d\ge 3$, then Amerik's result \cite[Corollary 2.2]{Am} implies that the Hodge numbers of $Y$ coincide either with that of
$X$ or with that of $\bP^n$. This excludes some rational homogeneous spaces $G/P$ as $Y$ with $\rho(G/P)=1$ \cite[Proposition 5.1]{BGG}. 
Again, if $d\ge 3$, then $Y$ can not be a hypersurface of degree~$\ge 3$ in $\bP^{n+1}$, by \cite[Corollary 1.1]{Choe}. Also if $d\ge n+3$ then there exists a constant $C(X)$ depending on discrete invariants of $X$ such that $Y$ can not be a smooth $n$-quadric hypersurface when $\deg f\ge C(X)$ \cite[Page~264]{Am}. 
\end{remark}

\section{Bounds of degrees of morphisms and conditions to be the projective space}

The following is essentially from \cite{HP} and needed for the proof of Proposition \ref{prop: Weil gen}. Here we have weakened the seemingly needed ampleness of the divisor $D$ omitted there, and used Khovanskii-Teissier equality case and Fujita's $\Delta$-genus theory for singular varieties to remove the klt assumption on $X$ there and slightly simplified the argument.

For a Weil divisor $D$ on a normal variety $X$, we use $\cO_X(D)$ to denote the corresponding reflexive sheaf; see \cite[Proposition~(2) and Theorem~(3), page 282]{Reid} for its properties. By abuse of notation, sometimes we use $H^0(X, D)$ to denote $H^0(X, \cO_X(D))$.

\begin{proposition}\label{Prop: HP} (cf. \cite[Theorem 3.2, Corollary 3.4]{HP})
Let $X$ be a normal projective variety of dimension $n$ and $D$ an integral $\bQ$-Cartier nef and big divisor such that $\dim H^0(X, \cO_X(D)) \ge n + D^n$, and the rational map
$\Phi_{|D|}$ is generically finite. Then:
\begin{enumerate}
\item
Fujita's $\Delta$-genus 
$\Delta(X, D) := n + D^n - \dim H^0(X, D)$ vanishes.
There are birational morphisms $\mu: V \to X$ and
$h: V \to W$ between normal projective varieties with $V$ being smooth, and very ample divisor $N$ on $W$, such that
$\mu^*D = h^*N$, so $D$ is Cartier and $|D|$ is base point free.
\item
If $D$ is ample, then $h$ can be identified with $\mu$ (and $W = X$).
\item
If $D$ is ample and $D^n \le 1$, then 
$X = \bP^n$ and $\cO_X(D) = \cO(1)$.
\end{enumerate}
\end{proposition}

\begin{proof}
By Hironaka's resolution theorem, there is a birational morphism $\mu: V \to X$ such that $V$ is smooth projective, $$\mu^*D = M + E$$
where $|M|$ (the moving part) is base point free with $H^0(V, M) \cong H^0(X, D)$, and $E$ is effective.
Since $|D|$ gives a generically finite map by the assumption, the base point free $|M|$ gives a generically finite morphism; thus $M$ is big and nef.
By the nefness of $\mu^*D$ and $M$ and the effectivity of $E$, we have
\begin{equation}\tag{*}\label{TE}
D^n = (\mu^*D)^n \ge (\mu^*D)^{n-1} M \ge \cdots \ge M^n
\end{equation}
Now the non-negativity of Fujita's $\Delta$-genus (see \cite[Theorem 1.1]{Fuj89}) and our assumption imply
$$0 \le \Delta(V, M) = n + M^n - \dim H^0(V, M) = n + M^n - \dim H^0(X, D)  \le n + D^n - \dim H^0(X, D) \le 0.$$
Thus the above and all the ones in the \eqref{TE} are equalities. Hence $\Delta(V, M) = 0 = \Delta(X, D)$, and
$$\begin{aligned}
(\mu^*D)^{n-1} . M &= M^n = D^n,  \\
((\mu^*D)^{n-1}M)^n &= (\mu^*D^n)^{(n-1)} M^n.
\end{aligned}$$
Now the last equality case of Kovanskii-Tessier implies that $M$ and $\mu^*D$ (and hence also $E$) are proportional in 
the $\bR$-coefficient N\'eron-Severi space (cf. \cite[Proposition~3.5]{LX}).

If $E \ne 0$ then, being proportional to $M$, is also nef and big, so the first inequality in the \eqref{TE} above is strict, which is a contradiction.
Thus $E = 0$. Hence $\mu^*D = M$ is Cartier, so is $D = \mu_*M$. On the other hand, by the vanishing of $\Delta(V, M)$, we have $$M = h^*N$$ for some birational morphism $h: V \to W$ to a normal projective variety and a very ample divisor $N$ (cf.~\cite[Theorem 1.1]{Fuj89}). 
This proves the assertion (1). 

The assertion (2) is clear by the rigidity lemma and ampleness of $D$ and $N$
(cf. \cite[Lemma 1.15]{Deb}).
For (3), (1) and (2) and the push-forward by $h = \mu$ imply (on $V = X$)
$$D = N,$$
so $D$ is integral, Cartier and ample. Now the assumption $D^n \le 1$ means $D^n = 1$.
Hence $(X, \cO(D)) = (\bP^n, \cO(1))$ by \cite[Corollary 4.3]{Fuj82} or \cite[Theorem~1.1]{KO}.
\end{proof}

Below is another key proposition which is needed for Theorems \ref{thm: e lc'}, \ref{thm: e lc-1}, and \ref{thm: e lc-2}.

\begin{proposition}\label{prop: Weil gen}
Let $X \subset \bP^{n+1}$ be a very general hypersurface of degree $d$ and $f: X \to Y$ a surjective morphism of degree $m \ge 2$ to a normal projective variety $Y$.
Assume either $n \ge 3$ and $d \ge 3$, or $n = 2$ and $d \ge 4$. Assume additionally $\pi_1(Y_{\rm reg})=(1)$ (which holds whenever $Y$ is smooth or $m := \deg f$ is prime by Proposition \ref{prop: basic}). Write $\Pic(X) = \bZ[H_X]$ with $H_X$ a general section of $\cO_{\bP^{n+1}}(1)|_X$ as in Proposition \ref{prop: basic}.
Let $L_Y$ be an integral divisor on $Y$ such that 
$$f^*L_Y \sim t_y H_X$$
for some minimal integer $t_y \ge 1$ (so that $L_Y$ is ample). Then:
\begin{enumerate}

\item 
For every integral Weil divisor $M$ on $Y$ we have $M \sim sL_Y$ for some integer $s$.

\item
Suppose $Y$ has only $\varepsilon$-lc singularities and $-K_Y$ is ample.
If $t_y = 1$, i.e., $f^*L_Y \sim H_X$, then $\deg f \le D_1 d$ for some constant $D_1 = D_1(n, \varepsilon)$ (with $D_1 = 1$ when $L_Y^n$ is an integer).

\item
Suppose that
 $S(f) :=\{ H\in |H_X|\, ;\,  \deg f|_{H}= \deg f \}$ is non-empty so that $f_*H = mH_Y$ where $H_Y = f(H)$ for some $H \in S(f)$. Then $H_Y \sim L_Y$, $t_y = 1$ and $S(f) = f^*|L_Y|$ (a linear subspace of $|H_X|$).

 \item 
$S(f) \ne \emptyset$ holds if and only if
$t_y = 1$ and $|L_Y| \not= \emptyset$.

\item
Suppose that $m = \deg f \ge d$. Then ``\,${\rm Codim} \, S(f) \le 1$ in $|H_X|$'' holds if and only if $\deg f = d$, $Y = \bP^n$ and $\cO(L_Y) = \cO(1)$.
\end{enumerate}
\end{proposition}

\begin{proof}
(1) Write $f^*M \sim s H_X$ with integer $s \ne 0$, say $s \ge 1$ after replacing $M$ by $-M$ if needed. By the minimality of $L_Y$, we have $t_y | s$. So $f^*M \sim f^*(s/t_y)L_Y$. Using the projection formula to push forward by $f_*$, our $M$ is $\bQ$-linearly equivalent and hence linearly equivalent to $(s/t_y)L_Y$ by the assumption $\pi_1(Y_{\reg}) = (1)$ (see Proposition \ref{prop: basic} (10)).

(2) By Birkar \cite[Theorem~1.1]{Bir}, there is a constant $D_1' = D_1'(n, \varepsilon)$ such that $|D_1' L_Y|$ gives a birational map, so $(D_1' L_Y)^n \ge 1$. Taking top intersection, we get $d = H_X^n = (f^*L_Y)^n = (\deg f) L_Y^n \ge (\deg f)/D_1$, where $D_1 = D_1(n, \varepsilon) := (D_1')^n$.

(3) By Proposition \ref{prop: basic}, $\Pic(X) = \bZ[H] = \bZ[H_X]$, and every member of $|H_X|$ is reduced and irreducible.
By the projection formula and the assumption, we have $f^*H_Y = H \sim H_X$. Thus $H_Y \sim L_Y$ and $t_y = 1$ (cf. (1)). This shows $S(f) \subseteq f^*|L_Y|$. Conversely, suppose $H' = f^*H_Y'$ ($\sim f^*L_Y \sim f^*H_Y = H \sim H_X$) for some $H_Y' \sim L_Y$. 
By the projection formula, $f_*H' = m H_Y'$, and hence $H_Y'$ is reduced (and irreducible) since $m = \deg f$. So $H_Y' = f(H')$, and hence $\deg f|_{H'} = m = \deg f$. Namely, $H' \in S(f)$.

(4) The ``only if'' part follows from (3). Conversely, suppose $t_y = 1$ and $L_Y \ge 0$. Thus $H_{X}' := f^*L_Y$ is in $|H_X|$ and hence it and then also $L_Y$ are reduced and irreducible, so $L_Y = f(H_X')$. By the projection formula, $mL_Y = f_*H_X'$, so $H_X' \in S(f)$.

(5) Assume first ${\rm Codim} \, S(f) \le 1$ in $|H_X|$. By (3), we have $f^*L_Y \sim H_X$. Taking top intersection and by assumption, $d = H_X^n = (f^*L_Y)^n = m L_Y^n \ge d L_Y^n$. Thus $L_Y^n \le 1$. Since $S(f) = f^*|L_Y|$ (cf. (3)), we have $\dim |L_Y|  = \dim S(f) \ge \dim |H_X| -1 = n$. Thus
\begin{equation}\tag{*}\label{eq: degree2}
\dim H^0(Y, L_Y) \ge n + 1 \ge n + L_Y^n.
\end{equation}

\begin{claim}\label{claim: no fib}
$|L_Y|$, or equivalently $S(f)$, defines a generically finite map $\Phi_{|L_Y|}$ (resp. $\Phi_S$).
\end{claim}

If $S(f) = |H_X|$, then as $H_X$ is very ample, the claim is true. If $\dim S(f) = \dim |H_X| - 1$, the linearity of $S(f) = f^*|L_Y|$ and Proposition \ref{prop: basic}~(2) imply that $S(f)$ is the restriction to $X$ of a sub-linear system of $|\mathcal{O}_{\bP^{n+1}}(1)|$ defined by one linear equation in the dual space $(\bP^{n+1})^*$. Thus $S(f) = \{ H_P \in |\cO_{\bP^{n+1}}(1)| ~;~ Q\in H_P\}|_X$
for some point $Q \in \bP^{n+1}$. Hence the map $\Phi_S$ is the projection $\varphi: \bP^{n+1} \setminus \{Q\} \to \bP^n$ restricted to $X$. If $\varphi(X) = \bP^n$ then the claim follows. Otherwise, $\varphi(X)$ is a hypersurface in $\bP^n$ in coordinates $X_0, \dots, X_n$ (while $\bP^{n+1} \supset X$ is in coordinates $X_0, \dots, X_{n+1}$) so that $X$ is the cone over $\varphi(X)$ and hence defined as $X = V(F(X_0, \dots, X_n))$ for some homogeneous polynomial $F$ in the first $n+1$ coordinates of $\bP^{n+1}$. This contradicts $X$ being a general hypersurface. The claim is proved.

We return back to the proof of the proposition.
With Claim \ref{claim: no fib} and the $(*)$ above, Proposition \ref{Prop: HP} implies $Y = \bP^n$ and $\cO(L_Y) = \cO(1)$. Also $d = H_X^n = (f^*L_Y)^n = \deg f$.

Conversely, if $m = d$, $Y = \bP^n$ and $\cO(L_Y) = \cO(1)$, 
then $t_y^n d = (t_y H_X)^n = (f^*L_Y)^n = m$, so $t_y = 1$.
Thus, by (4) (3), $S(f) = f^*|L_Y|$ and
$\dim S(f) = \dim H^0(Y, L_Y) - 1 = n = \dim |H_X| - 1$. 

This proves the proposition.
\end{proof}

\begin{remark}
In the case of Proposition \ref{prop: Weil gen} (5), if we assume $d\ge 2n+3$ then Lemma~\ref{lem: degree} 
(whose assumption is satisfied because $\deg f = d$ as in the proof of Propostion~\ref{prop: Weil gen}) shows that the fiber $Z_y$ for a general point $y\in Y$ of $f$ lies on the unique line $\ell_y$ in $\bP^{n+1}$. In fact \cite[Lemma 2.4]{BCD} says that $f^{-1}(y)$ lies on the unique line $\ell_y$ if all points in $f^{-1}(y)$ are distinct points. Then imitating \cite[Proof of Lemma 4.1]{BCD} and \cite[Proof of Theorem C]{BDet}, 
one might be able to show that the morphism $f: X \to Y$ is nothing but a projection from a point in $\bP^{n+1} \setminus X$. However, we do not pursue a proof here, since this paper focuses on the characterization of $Y$.
\end{remark}

When $\deg f = p$ is prime, there is an upper bound for $\deg f$ which depends on $\deg X = d$ linearly.

\begin{theorem}\label{thm: e lc'}
Let $X \subset \bP^{n+1}$ be a very general hypersurface of degree $d$ 
and $f: X \to Y$ a surjective morphism of prime degree to a normal projective variety $Y$
with at worst $\varepsilon$-lc singularities. Assume either $n \ge 3$ and $d \ge 3$, or $n = 2$ and $d \ge 4$. Assume further either $Y$ is factorial (or just $L_Y^n$ is an integer) or $-K_Y$ is ample, when $f^*L_Y \sim H_X$. Then:
\begin{enumerate}
\item 
There are constants $D_j = D_j(n, \varepsilon)$ with $D_1 \ge 1$, with $D_1 = 1$ when $Y$ is factorial (or just $L_Y^n$ is an integer)
and with $D_2 = (D_2')^{1/(n-1)}$ in terms of upper bound $D_2' = D_2'(n, \varepsilon)$ of $(-K_Y)^n$
(as in Birkar~\cite[Theorem 1.1]{Bir}, or \cite[Corollary~1]{KMM} with $D_2' = \{n(n+1)\}^n$ when $Y$ is smooth), 
such that 
$$\deg f \, \le \, {\rm max} \,\{D_1 d, \,\,\,\,\frac{D_2}{d^{1/(n-1)}}\}.$$
\item 
If $Y$ is factorial (or just $L_Y^n$ is an integer) and $d \ge D_2^{(n-1)/n}$, then $\deg f \le d$.
\end{enumerate}
\end{theorem}

\begin{proof}
Note that (2) follows from (1). We prove (1).
Let $p = \deg f$ which is assumed to be a prime number. We may assume $p \ge d$. By our assumption, all conditions in Proposition~\ref{prop: basic} 
are satisfied. So $Y$ is $\bQ$-factorial and klt with $\rho(Y) = 1$, and $\Pic(X) = \bZ[H_X]$ with $H_X$ a general section of $\cO_{\bP^{n+1}}(1)|_X$; every member of $|H_X|$ is reduced and irreducible.

Consider $S(f) :=\{ H\in |H_X|\, ;\,  \deg f|_{H}= \deg f \}$. If $\dim S(f) \ge n$, then $p = d$, $Y = \bP^n$ and $\cO(L_Y) = \cO(1)$, by Proposition \ref{prop: Weil gen}. Thus we may assume $\dim S(f) \le n-1$, i.e., $S(f)$ is of codimension $\ge 2$ in $|H_X|$.

Let $H_Y := f(H_X)$, an integral and $\bQ$-Cartier ample divisor.
Let $f_{H_X} = f|_{H_X} : H_X \to H_Y$ and $m_{H_X} = \deg f_{H_X}$.
Write $f^*H_Y \sim b_{H_X} H_X$. By the projection formula,
$p H_Y = f_*f^*H_Y \sim b_{H_X} m_{H_X} H_Y$.
So $p = b_{H_X} m_{H_X}$. Since $p$ is prime, $(m_{H_X}, b_{H_X}) = (p, 1)$ or $(1, p)$. Since we assume $S(f)$ has codimension $ \ge 2$ in $|H_X|$, we have $(m_{H_X}, b_{H_X}) = (1, p)$ for general $H_X$.
Namely, $f_*H_X = H_Y$ and $f^*H_Y \sim pH_X$. Write $H_Y \sim h_yL_Y$. Then $h_yt_y H_X \sim f^*(h_yL_Y) = f^*H_Y \sim pH_X$. Hence $h_yt_y = p$. So $(h_y, t_y) = (p, 1)$ or $(1, p)$.

Consider the case $(h_y, t_y) = (p, 1)$. Then $H_Y \sim pL_Y$ and $f^*L_Y \sim H_X$. If $L_Y^n$ is an integer (so $\ge 1$, since $L_Y$ is ample), we have
$$d = H_X^n = (f^*L_Y)^n = p L_Y^n \ge p.$$ If $-K_Y$ is ample, by Proposition \ref{prop: Weil gen} (2), we have $$p \le D_1 d$$ for some constant $D_1 = D_1(n, \varepsilon)$.

Consider the case $(h_y, t_y) = (1, p)$. Then $H_Y \sim L_Y$ and $f^*L_Y \sim pH_X$.

\begin{claim}
$-K_Y$ is ample.
\end{claim}

We prove the claim. Suppose the contrary that $-K_Y$ is not ample. Then by Theorem \ref{thm: main2}, $K_Y$ is ample. By the ramification divisor formula, $(d-n-2)H_X \sim K_X = f^*K_Y + R_f$. Pushing forward using $f_*$, we get
$(d-n-2)H_Y \ge pK_Y \ge pL_Y \sim pH_Y$. So $p \le d-n-2$, contradicting the assumption  $p \ge d$. This proves the claim.

\medskip

We return back to the proof of the theorem. Write $-K_Y \sim rL_Y \sim rH_Y$ for some integer $r \ge 1$.
Then the ampleness of them, imply 
\begin{equation}\tag{*}\label{eq: deg bound}
p^nd = (pH_X)^n = (f^*H_Y)^n = p H_Y^n \le p(-K_Y)^n \le p D_2'(n, \varepsilon)
\end{equation}
for some constant $D_2' = D_2'(n, \varepsilon)$ as specified in the statement of the theorem.
Now (1) follows from the $(*)$ above by letting
$D_2 = (D_2')^{1/(n-1)}$.
\end{proof}


\begin{proof}[{\bf Proof of Theorem~\ref{thm: e lc-1} and Remark \ref{rem: on main thm} (4)}]

Theorem \ref{thm: e lc-1} (1) follows Proposition \ref{prop: Weil gen} (4) and (2). For Theorem \ref{thm: e lc-1} (2), by Theorem \ref{thm: main2}, we may assume that $\deg f \ge C(n, d)$ and hence $-K_Y$ is ample. Then it follows by taking the current $D_2$ to be the maximum of its name sake in Theorem \ref{thm: e lc'} and $C(n, d)$. For Theorem \ref{thm: e lc-1} (3) and Remark \ref{rem: on main thm} (4), we may assume $\deg f \ge d$, and then they follow from Proposition \ref{prop: Weil gen} (5).
\end{proof}

\begin{proof}[{\bf Proof of Theorem~\ref{thm: e lc-2} and Remark \ref{rem: on main thm} (5)}]
We assume either $n \ge 3$ and $d \ge 3$, or $n = 2$ and $d \ge 4$ (instead of $n \ge 3$ and $d \ge 3$).
Theorem \ref{thm: e lc-2} (1) follows Proposition \ref{prop: Weil gen} (4) and (2). Theorem \ref{thm: e lc-2} (2) implies (3), while (2) itself follows from Theorem \ref{thm: e lc'} with $\varepsilon = 1$ since $Y$ is klt and factorial and hence has at worst canonical singularities.
\end{proof}

As consequences of Theorem \ref{thm: e lc'}, we now give sufficient conditions for the target $Y$ of a very general hypersurface $X \subset \bP^{n+1}$ to be $\bP^n$. We start with dimension $3$ and then dimension $\ge 4$.

Let $Y$ be a smooth Fano $n$-fold with $\Pic Y = \bZ[H]$. Write $-K_Y \sim rH$. We call $r$ the {\it Fano index} of $Y$ and $n+1-r$ the {\it Fano co-index} of $Y$. It is known that $r \le n+1$, and $r=n+1$ (resp. $r = n$) holds if and only if $Y \cong \bP^n$ (resp. $Y$ is a hyperquadric in $\bP^{n+1})$; see \cite[Corollary to Theorem 1.1, Corollary to Theorem 2.1]{KO}.

\begin{corollary}\label{cor: p3}
In Theorem~\ref{thm: e lc-2}, assume further
\begin{enumerate}
\item $\deg f$ is prime, 
\item $Y$ is smooth,
\item $\deg f \ge d$,
\item $n=3$, and
\item $d > n(n+1) = 12$.
\end{enumerate}
Then $\deg f = d$ and $Y \cong \bP^3$.
\end{corollary}

\begin{proof} 
Since $Y$ is smooth, $\Pic Y = \bZ[L_Y]$ by Propositions \ref{prop: Weil gen} and \ref{prop: basic}. Let $p = \deg f$, a prime number. 
If $\dim S(f) \ge n$, then $p = d$, $Y = \bP^n$ and $\cO(L_Y) = \cO(1)$, by Proposition \ref{prop: Weil gen}. Thus, we may assume $\dim S(f) \le n-1 = 2$.
By Theorem \ref{thm: e lc-2} and Remark~\ref{rem: on main thm} (2), we have $\deg f \le d$, so $\deg f = p = d$. Since $\dim S(f) \le n-1$, by the proof of Theorem \ref{thm: e lc'}, we may assume that $(m_{H_X}, b_{H_X}) = (1, p)$ for general $H_X$, i.e., $f_*H_X = H_Y := f(H_X)$ and $f^*H_Y \sim pH_X$; further, $(h_y, t_y) = (p, 1)$ or $(1, p)$, i.e., $H_Y \sim pL_Y$ and $f^*L_Y \sim H_X$, or $H_Y \sim L_Y$ and $f^*L_Y \sim pH_X$.

By Theorem \ref{thm: main2}, $K_Y$ or $-K_Y$ is ample. Thus $-K_Y$ is ample (cf.~\cite[Proof of Lemma 3.0]{Am}).

Suppose $H_Y \sim L_Y$ and $f^*L_Y \sim pH_X$. Then, $p^{n+1} = p^n d = (f^*L_Y)^n = p L_Y^n$. Also, by \cite[Corollary 1]{KMM}, we have (with $n = 3$)
$$(n(n+1))^n \ge (-K_Y)^n\ge L_Y^n=p^n = d^n$$
contradicting our assumption. 

Thus, $H_Y \sim pL_Y$ and $f^*L_Y \sim H_X$,
hence $L_Y^3=1$. Write $-K_Y \sim rL_Y$. We may assume $r \le n = 3$. Then $(-K_Y)^3 = r^3$.
Using the classification theory of Fano threefolds with Picard number $\rho=1$ (cf. \cite[Table 12.2]{Par}), we have $r = 2$ and $\dim H^{1,2}(Y) = 21$ (so $Y$ is a del Pezzo threefold of degree 1). Note that the Hodge pair $(h^{3,0}, h^{1, 2})$ of $Y$ is $(0, 21)$ while that of $X$ is $(+, *)$ (because $d > n+1$) and of $\bP^3$ is $(0, 0)$. We get a contradiction to \cite[Corollary 2.2]{Am}.
\end{proof}

\begin{corollary}\label{cor: pn}
In Theorem~\ref{thm: e lc-2}, assume further
\begin{enumerate}
\item $\deg f$ is prime, 
\item $Y$ is a smooth Fano variety,
\item $\deg f \ge d$,
\item $n\ge 4$,
\item $d > n(n+1)$, and 
\item Fano index of $Y$ is $\ge n-2$ and $\ne n-1$.
\end{enumerate}
Then $\deg f = d$ and $Y \cong \bP^n$.
\end{corollary}

\begin{proof} 
Since $Y$ is smooth, $\Pic Y = \bZ[L_Y]$ by Propositions \ref{prop: Weil gen} and \ref{prop: basic}. Let $p = \deg f$, a prime number. 
If $\dim S(f) \ge n$, then $p = d$, $Y = \bP^n$ and $\cO(L_Y) = \cO(1)$, by Proposition \ref{prop: Weil gen}. Thus we may assume $\dim S(f) \le n-1$. By the same argument as in Corollary \ref{cor: p3}, we may assume $p = d$, $H_Y \sim pL_Y$, $f^*L_Y \sim H_X$ (so $t_y = 1$), and $L_Y^n = 1$. 
On the other hand, by Ambro \cite[Main Theorem]{Ambro}, we have $\dim |L_Y|\ge n-1$. Now $S(f) = f^*|L_Y|$ by Proposition \ref{prop: Weil gen}. Therefore, $\dim S(f) = \dim |L_Y| = n-1$. So Fujita's $\Delta$-genus of $(Y, L_Y)$ is $1$, which would reach a contradiction.

Indeed, by applying Fujita \cite[Theorem 2.1 and Proposition 3.2]{Fuj82}, a general member $D$ in $|L_Y|$ is irreducible and smooth. The polarized pair $(Y, L_Y)$ satisfies the following:
\begin{enumerate}
\item $Y$ is a smooth Fano variety with Fano index $\ge n-2$ and $\ne n-1$, i.e., of Fano co-index $\le 3$ but $\not= 2$, and with $\Pic Y = \bZ[L_Y]$,
\item $L_Y^n=1$,
\item $\dim |L_Y|=n-1$, 
\item Fujita's $\Delta$-genus $\Delta(Y, L_Y)=1$, and
\item $\Pic D=1$ for a general member $D$ in $|L_Y|$ (cf. \cite[Example 3.1.25]{Laz}).
\end{enumerate}

Then the pair $(D, L_Y|_D)$ satisfies again all the above conditions (1)--(5). Indeed,  by the adjunction formula, $D$ is a smooth Fano variety of dimension $n-1$ and Fano index $\ge (n-1)-2$. Since $Y$ is a smooth Fano variety, $H^1(Y, \cO_Y)=0$ by Kodaira vanishing. Therefore, we keep
$$\Delta(D, L_Y|_D) = \Delta(Y, L_Y)=1$$
by \cite[(1.5)]{Fuj82}. And by \cite[Example 3.1.25]{Laz}, the restriction $\Pic(Y)\to \Pic(D)$ is an isomorphism because $n\ge 4$.  Therefore $L_Y|_D$ is an ample generator of $\Pic D$ and $(L_Y|_D)^{n-1}=1$. Thus $Y$ and $D$ share the same Fano co-index.
Also from the exact sequence
$$0\to \cO_Y\to \cO_Y(L_Y)\to \cO_D(L_Y|_D)\to 0,$$
we have $\dim |L_Y|_D|=\dim |L_Y|-1$ because $H^1(Y, \cO_Y)=0$. It gives $\dim |L_Y|_D|=n-2$. Then we can apply again \cite[Theorem 2.1 and Proposition 3.2]{Fuj82}.

This process reduces to the case where $\dim Y = 4$. We note that this process keeps the Fano co-index. Then we have a smooth Fano threefold $D\in |L_Y|$ again by \cite[Theorem 2.1 and Proposition 3.2]{Fuj82}. This $D$ satisfies the above conditions (1)--(4). As in Corollary \ref{cor: p3}, write $-K_D \sim r L_Y|_D$ with $r \le \dim D + 1 = 4$, and use $(-K_D)^3 = r^3$ and \cite[Table 12.2]{Par}. We see that either $r = 4$ and $D = \bP^3$, or $r = 2$. If $D = \bP^3$, then ($1  = $) $\Delta(D, L_Y|_D)=\Delta(\bP^3, \cO(1))=0$, a contradiction. If $r = 2$, then $D$ and hence $Y$ have Fano co-index 2, which we exclude by the assumption. This proves the corollary.
\end{proof}

\section
{More conjectures and relations among them}

In this section we will discuss a related
conjecture about dominant rational maps. Precisely, one may also consider the following one
about dominant rational maps from a very general hypersurface (cf. Guerra and Pirola in \cite[Conjecture, page 290]{gp}):

\begin{conjecture}\label{conj: 1}
Let $X\subset \bP^{n+1} $ be a very general hypersurface of degree $d\ge n+3$, and $Y$ a smooth projective $n$-fold. Suppose there is a dominant rational map $f: X\dashrightarrow Y$ with ${\rm deg}\, f>1$. Then $Y$ is rationally connected.
\end{conjecture}

By Hodge theory (cf.~Guerra and Pirola \cite[Proposition 3.5.2]{gp}), to prove Conjecture \ref{conj: 1}, it suffices to consider dominant rational maps $f: X\dashrightarrow \widetilde{Y}$,
where $\widetilde{Y}$ is a smooth model of $Y$, which is simply connected and without holomorphic global $n$-forms, i.e., $p_g(\widetilde{Y})=0$.

In \cite{gp}, Guerra and Pirola proved that the $Y$ in Conjecture \ref{conj: 1} is not of general type for the case $n=2$ and $5\le d\le 11$. 
Conjecture \ref{conj: 1} for $n=2$ is proved in the joint paper \cite{LP15} of Pirola and the first-named author of the present paper.
In Conjecture \ref{conj: 1}, the assumption `general' is clearly necessary. Indeed, if one chooses a special surface $X$ in $\bP^3$, it might admit a dominant rational map to a surface of general type $Y$. Classical Godeaux surfaces, which are of general type, are obtained as $\bZ/(5)$-quotients of $\bZ/(5)$-invariant quintics (cf.~Godeaux \cite{Go}). Also for instance, a Fermat hypersurface of degree $d$ in $\PP^{n+1}$ with $d=d_1a$ can map to a Fermat hypersurface of degree $d_1$ in $\PP^{n+1}$, which is not rationally connected when $d_1 \ge n+2$, by taking $a$-th power.

In \cite{LP15}, the assumption `very general' is used to deal with the case $n=2$ and to exclude countably many different deformation types of non-ruled simply connected minimal surfaces with $p_g=0$.

Here, a projective variety is {\it rationally connected} in the sense of Campana and Koll\'ar-Miyaoka-Mori,
if two general points are connected by a rational curve; {\it this is a birational notion} (cf.~Koll\'ar \cite[VI, Proposition 3.3]{Ko}). 
We consider `rational connectedness' because this assumption fits well with the birationality and induction argument on the dimension. 

We remark that one approach to prove Conjecture~\ref{conj: 1} is by induction via the following conjecture 
and Proposition~\ref{prop: surjective}. Indeed, Conjecture~\ref{conj: 1} in dimension $\le n-1$ and Conjecture~\ref{conj: induction} in dimension $n$ imply Conjecture~\ref{conj: 1} in dimension $n$.

\begin{conjecture}\label{conj: induction}
Let $X\subset \bP^{n+1} $ be a very general hypersurface of degree $d\ge n+3$, and $W$ a smooth projective $n$-fold which is non-uniruled.
Let $H_X := X\cap \, H$, where $H$ is a very general hyperplane
in $\bP^{n+1}$.
Suppose there is a dominant rational map $h: X \dashrightarrow W$ with ${\rm deg}\, h>1$.
Then the restriction of $h$ to $H_X$ is a non-birational map from $H_X$ to its image in $W$.
\end{conjecture}

We remark that although there are countably many different deformation types for simply connected minimal non-ruled surfaces with $p_g=0$, the dimension of each moduli space of simply connected minimal non-ruled surfaces with $p_g=0$ is bounded by 19 (cf. \cite[Corollary~2.5.3]{gp} for surfaces of general type, \cite[Proof of Theorem 1.1]{LP15} for surfaces with Kodaira dimension $\kappa=1$). This helps to solve the case $n=2$ for Conjectures~\ref{conj: 1} and \ref{conj: 2} in \cite{LP15}.
However, when $\dim Y\ge 3$, there are no known results for the dimension of moduli of simply connected minimal non-ruled varieties with $h^i(Y, \cO_Y)=0$ for all $i >0$.

\begin{proposition}\label{prop: surjective}
Let $g: X \dasharrow Z$ be a dominant rational map between projective varieties over an algebraically closed field of characteristic zero, where $X$ is normal and $\dim\, X > \dim\, Z\ge 1$. Let $H$ be a nef and big Cartier divisor, and $\Lambda$ a non-trivial linear subsystem of $|H|$ with ${\rm Bs}\,\Lambda = \emptyset$.
Then for general $H_1, \dots, H_r$ in $\Lambda$ with $r = \dim\, X - \dim\, Z$, the map $g$ restricts to the intersection $H_1 \cap \ldots \cap H_r$ as a generically finite dominant map.
\end{proposition}

\begin{proof}
Replacing $Z$ and $X$ by smooth models $Z'$ and $X'$, and pulling back $H$ and $\Lambda$ to $X'$, we may assume that $X$ and $Z$ are smooth and $g$ is a surjective morphism. 

Since the Iitaka dimension $\kappa(X, H) = \dim\, X > \dim\, Z$, any $H'$ in $\Lambda$ dominates $Z$, or else $H' \le f^*D$ for some divisor $D$ on $Z$ and 
$$\kappa(X, H) = \kappa(X, H') \le \kappa(X, f^*D) = \kappa(Z, D) \le \dim\, Z,$$
which is a contradiction.
By Bertini's theorem, a general member $H_1$ of $\Lambda$ is irreducible and smooth as $\Lambda$ is base point free and $\dim X \ge 2$. Replace $X$ by $H_1$ and $\Lambda$ by $\Lambda \, |_{H_1}$, which remains a non-trivial base point free linear subsystem of $|H_1 |_H|$ with $H_1|_H$ being nef and big. By induction, we have a surjective morphism 
$$H_1 \cap \ldots \cap H_r \to Z$$ 
between varieties of the same dimension, thereby proving the proposition.
\end{proof}

To end this section, we remark that there is a lower bound for $\mathrm{deg}~f$ if the generically finite map $f: X\dashrightarrow Y$ is a non-birational dominant map.
Lemma~\ref{lem: degree} below is essentially from \cite{BCD}.

\begin{lemma} \label{lem: degree}
Let $X\subset \bP^{n+1} $ be a very general hypersurface of degree $d\ge n+3$ with $n\ge 2$, and $Y$ a smooth projective $n$-fold. Suppose there is a dominant rational map $f: X\dashrightarrow Y$ with $\mathrm{deg}~f >1$.
Then the following assertions hold.
\begin{itemize}
\item[(1)] $\mathrm{deg}~f \ge d-n$.
\item[(2)] If $\mathrm{deg}~f\le 2d-2n-3$, then the fiber $Z_y$ for a general point $y\in Y$ of $f$ lies on a line $\ell_y$ in $\bP^{n+1}$ and $\mathrm{deg}~f\le d$.
\end{itemize}
\end{lemma}

\begin{proof}
By Claim 2.4,
$p_g(Y)=0$. Therefore, given $\omega \in H^0(X, \cO(K_X))$, its \emph{trace} (cf. \cite[Definition~2.1]{BCD}) in $H^0(Y, \cO(K_Y))$ vanishes. Accordingly, the general fiber of $f$ has the Cayley-Bacharach property 
(cf.~Bastianelli-Cortini-De Poi \cite[Definition~2.2]{BCD}) with respect to $\cO(K_X) = \cO_X(d-n-2)$. Namely, if a section of $\cO(K_X)$ vanishes at all but one point of the general fiber then it vanishes at every point of the general fiber. Hence we have $\mathrm{deg}~f \geq d-n$; see \cite[Proposition 2.3, Lemma~2.4]{BCD} or \cite[Theorem 2.3]{BDet}. This proves (1).

For (2), suppose that $\mathrm{deg}~f\leq 2(d-n-2)+1$. The rational map $f$ gives a correspondence $Z\subseteq X\times Y$ with vanishing trace map from $H^0(X, \cO(K_X))$ to $H^0(Y, \cO(K_Y))$. By \cite[Lemma 2.4 or Theorem 2.5]{BCD}, general fiber $Z_y$ of $f$ lies on a line $\ell_y\subset \bP^{n+1}$. Since $\deg X=d$, we have $\mathrm{deg}~f\leq d$ because $Z_y$ lies on a line $\ell_y$. 
\end{proof}

\end{document}